\documentclass[v,journal]{IEEEtran}
\usepackage{amsmath,amsfonts}
\usepackage{algorithmic}
\usepackage{algorithm}
\usepackage{array}
\usepackage[caption=false,font=normalsize,labelfont=sf,textfont=sf]{subfig}
\usepackage{textcomp}
\usepackage{stfloats}
\usepackage{url}
\usepackage{verbatim}
\usepackage{graphicx}
\usepackage{cite}
\usepackage{enumitem}
\usepackage{xurl}
\usepackage{placeins}
\usepackage{xcolor}
\begin{document}

\title{Grid-Interactive Operation of Solar-Integrated Data Centers for Coordinated Local and System-Level Decarbonization}

\author{Mingxu~Yang, Weiqi~Zhang, Hui~Geng, and~Jiaze~Ma

\thanks{M. Yang and J. Ma are with the City University of Hong Kong, Tat Chee Avenue, Hong Kong (e-mail: mingxyang9-c@my.cityu.edu.hk; jiazema@cityu.edu.hk)\emph{(Corresponding author: Jiaze Ma.)}.}%
\thanks{W. Zhang is with Uber Technologies Inc., 1655 3rd St, San Francisco, CA 94158 USA (e-mail: weiqiz@uber.com).}%
\thanks{H. Geng is with the University of Wisconsin–Madison, 1415 Engineering Dr., Madison, WI 53706 USA (e-mail: hgeng7@wisc.edu).}}

\markboth{Journal of \LaTeX\ Class Files,~Vol.~14, No.~8, August~2021}%
{Shell \MakeLowercase{\textit{et al.}}: A Sample Article Using IEEEtran.cls for IEEE Journals}


\maketitle

\begin{abstract}
The exponential growth of AI is accelerating the deployment of data centers (DCs), placing unprecedented strain on power infrastructures. In response, major IT corporations are increasingly adopting on-site solar generation to reduce grid dependence and meet sustainability targets. However, the true impacts of this strategy remain ambiguous. While DCs are flexible assets capable of temporal load-shifting, anchoring them to self-generated power may inadvertently constrain their grid responsiveness. To evaluate these trade-offs, we propose a receding-horizon optimization (RHO) framework coordinating job scheduling, grid interactions, and on-site solar generation for a stand-alone DC. Our findings reveal a critical paradox: although solar integration increases energy self-sufficiency and reduces overall DC emissions, it inherently limits the facility's capacity to absorb low-cost, low-carbon electricity from the grid. This implies a fundamental tension between individual corporate sustainability goals and system-wide grid decarbonization.
\end{abstract}

\begin{IEEEkeywords}
Data center, on-site solar, energy market, green computing
\end{IEEEkeywords}

\section{Introduction}
\IEEEPARstart{D}{riven} by the rapid development of AI, the electricity consumption of DCs in the U.S. has reached 176 TWh in 2024, representing 4.4\% of total U.S. electricity consumption\cite{shehabi20242024}. This is equivalent to the power used by 6.7 million average American households. Meanwhile, it is projected that electricity consumption by DCs in the U.S. will increase to two to three times the 2024 level by 2028 driven by the explosion of generative AI\cite{shehabi20242024}. In 2024, the total greenhouse gas emissions (GHGs) from U.S. DCs reached 61 billion kg of CO$_{2}$, which is comparable to the annual GHGs of 11 million cars on the streets \cite{shehabi20242024, EPA_PassengerVehicle_2025}. The huge electricity demand from DCs exerts substantial pressure on regional power supply and grid stability\cite{McLaughlin_DataCenterBoom_2025}. Unfortunately, the expansion of power grid infrastructure cannot keep pace with the exponential growth of DCs, leading to delays in project deployment\cite{IEA_Electricity_2025,  EPRI_3002028905}.

To mitigate the impact of DCs' electricity demand on grid operations, electricity markets have introduced peak-related penalties. DC electricity bills typically consist of usage charges, which depend on the amount of electricity used and the hourly electricity prices, and peak demand charges, which depend on the largest electricity load of a DC over a period of time, such as the peak power load in a month \cite{liu2013data,liu2014pricing}. Using Google DCs as an example, the peak demand charges can account for up to 56\% of the monthly electricity bill \cite{dabbagh2017shaving}. As a result, large IT enterprises have now recognized that exclusive reliance on the public power grid is insufficient to ensure power supply reliability and business continuity. This trend leads to the adoption of self-supply and diversified power procurement strategies, such as deploying on-site solar power to mitigate the reliance on the public grid. In addition to economic and reliability considerations, the development of renewable assets improves companies' Environmental, Social, and Governance (ESG) performance, which in turn affects company image, investor confidence, and market valuation \cite{galema2025esg}. Motivated by cost-effectiveness, reliability, and ESG considerations, major IT companies are turning to their own renewable and low-carbon energy sources, such as solar and nuclear power. For instance, Amazon has matched its DC electricity consumption through more than 52 renewable energy projects across Europe \cite{Amazon_CarbonFreeEnergy_2025}. Microsoft and Google have secured electricity supply through multiple solar power agreements\cite{Leong_MicrosoftZelestraSolar_2025, Constellation_CraneCleanEnergy_2024, DeChant_GoogleDataCenterPower_2026}. In an extreme case, Google have even been involved in restarting a few nuclear power plants or have entered into agreements with existing nuclear energy providers to supply their DCs with low-carbon electricity \cite{Wehner_GoogleNuclearDeal_2025}. Overall, the rapid expansion of DCs in the U.S. is placing increasing pressure on traditional power grids, driving major IT enterprises to pursue on-site renewable generation to ensure power reliability, cost effectiveness, and ESG performance.

The on-site renewable generation can mitigate the reliance on the public power grid, but the inherent variability of renewable energy sources also exposes DCs to power fluctuations \cite{gnibga2023renewable, Bojek_SolarPVGlobalSupplyChains_2022}. In parallel, the volatility of electricity prices adds another layer of complexity to the operation of the DC. Taking California as an example, electricity prices exhibit large fluctuations, as the highest regional prices can reach up to 1000 USD/MWh, while prices in the lowest-priced regions may even fall below zero \cite{CAISO_TodaysOutlookPrices_2025}. DCs' jobs fluctuate across diverse service types, as the incoming computing jobs are random \cite{gnibga2023renewable}. The joint volatility of renewable generation, electricity prices, and computational demand therefore renders DC scheduling inherently complex. At the same time, DC job scheduling exhibits intrinsic flexibility, as a substantial share of jobs can be deferred or curtailed, enabling demand-side flexibility \cite{ghamkhari2013energy, zhang2023quantifying}. The integration of on-site renewable generation further complicates scheduling by requiring coordinated decisions between self-generation, grid procurement and DCs' job management. 

\begin{figure*}[!b]	\center{\includegraphics[width=0.8\textwidth]{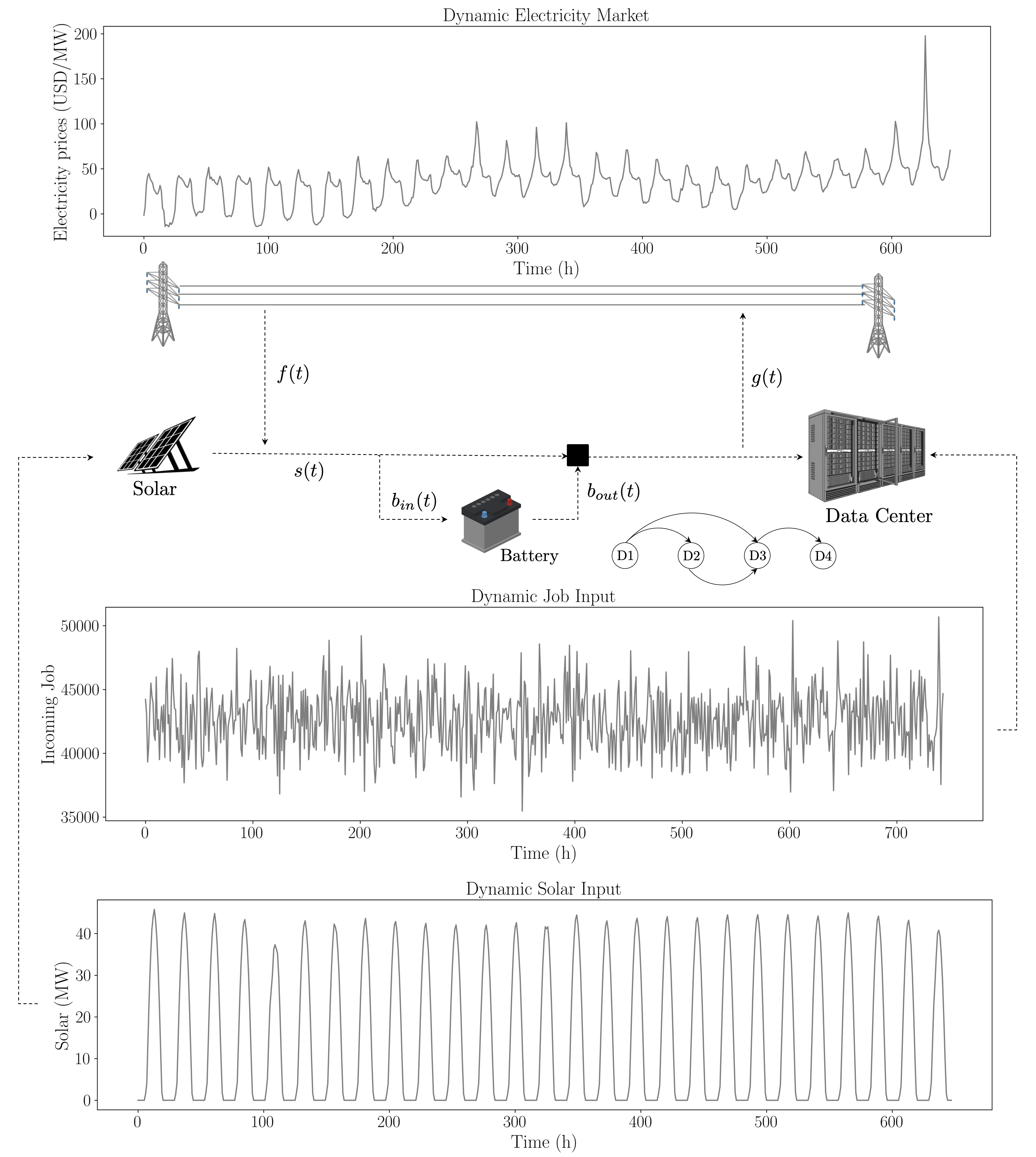}}
	\caption{DC scheduling under volatile electricity prices, solar inputs, and incoming computing jobs: DC has the flexibility to shift jobs based on solar availability and electricity prices with a battery buffer. Solar power can be selectively allocated to the DC, the power grid, or the battery. The DC is jointly powered by solar, battery, and the power grid.}
	\label{Schematic}
\end{figure*}

Motivated by the growing exposure of DCs to electricity markets and on-site renewable generation, a substantial body of literature has investigated job scheduling and energy management under volatile electricity prices, renewable output, and random job demand. Early studies first focused on improving quality of service (QoS) through job scheduling and server provisioning~\cite{beloglazov2012managing, delimitrou2013paragon, gu2020energy, sharma2020artificial}. As awareness of the massive energy footprint of DCs increased, the field advanced toward energy efficiency improvements that are sensitive to service quality~\cite{beloglazov2010energy, wang2015provision, yuan2021energy, ghamkhari2016energy, li2017energy}. Subsequent studies incorporated dynamic electricity prices and demand charges, demonstrating that load shifting and market-aware scheduling can substantially reduce operating costs~\cite{le2009cost, bahrami2018data, zhang2021hpc, zhang2024data}. Building on this line of work, our previous study~\cite{zhang2022exploring} proposed a RHO framework to coordinate DC job scheduling with dynamic electricity markets, capturing electricity price volatility and job flexibility in a rolling decision-making setting. Driven by stringent environmental policies and growing renewable energy penetration, recent studies have demonstrated that the integration of renewable energy sources can further reduce the DC's operating costs while strictly maintaining QoS guarantees~\cite{wang2014hierarchical, he2024analysis, lagana2018reducing, peng2022exploiting, agarwal2021redesigning}. However, how on-site renewable generation affects optimal job scheduling has received limited attention. DCs were also treated solely as electricity consumers, the impact of on-site renewable generation on the power grid is also unknown. As large-scale DCs increasingly deploy on-site solar generation, the interaction between market-driven scheduling incentives, non-dispatchable self-generated power, and prosumer operation introduces new operational trade-offs that are not fully captured by existing receding-horizon formulations. 

In this work, we extend our previous framework~\cite{zhang2022exploring} to jointly account for volatile electricity prices, dynamic solar output, and flexible job migration. Due to their inherent workload flexibility, DCs can act as temporally adaptable loads, shifting computational demand to absorb surplus, low-cost, and low-carbon electricity from the grid. However, the integration of on-site solar fundamentally reshapes this role by partially coupling energy consumption to local generation. While self-generation reduces reliance on grid electricity, it may also constrain the ability of DCs to respond to market signals, limiting their participation in low-price periods and weakening their effectiveness as demand-side flexibility resources. In addition, solar deployment introduces capital investment and embodied carbon emissions, further complicating its net environmental benefits. By anchoring operations to on-site generation, solar integration may inadvertently reduce the ability of DCs to utilize system-level low-carbon electricity that is abundant but temporally concentrated, such as during midday renewable peaks. As a result, it remains unclear whether on-site solar can simultaneously deliver economic savings and system-efficient decarbonization.

To address this gap, we develop a coordinated optimization framework that captures the interaction between flexible computing workloads, on-site solar generation, and electricity market dynamics under realistic operating conditions. Our results reveal a fundamental paradox: while integrating on-site solar reduces total greenhouse gas emissions and lowers electricity procurement costs at the facility level, it simultaneously increases the average carbon intensity of grid electricity consumption and diminishes the ability of data centers to act as flexible demand-side resources. This effect arises from the temporal alignment between local solar generation and grid-wide renewable availability, which limits the absorption of surplus low-carbon electricity. These findings highlight a misalignment between facility-level sustainability objectives and system-level grid decarbonization, underscoring the need for grid-interactive strategies that better coordinate distributed renewable generation with flexible demand.

\section{DC Scheduling and Energy Management Model}\label{sec:Methods}

We develop a computational framework for optimal DC job scheduling and energy management. Following the formulation in \cite{zhang2022exploring}, we model the system using a RHO framework. Upon receiving electricity prices, solar radiation forecasts, and incoming job information, the decision-making framework optimizes energy transactions, job execution, and energy flows over a 168-hour receding horizon. Only the first-step decision of each optimization horizon is implemented, while the remaining steps are updated in the next iteration. The final output is obtained by concatenating the first-step decisions across all horizons.

The DC consists of $\mathcal{I} = \{1, 2, \dots, I\}$ servers, and the scheduling horizon is discretized into time steps with a total length of $\mathcal{T} = \{1, 2, \dots, T\}$. Due to limitations in power supply, the operator does not activate all $|\mathcal{I}|$  servers at every time step. We define $I_t \subseteq I$ as the number of servers available at time t. Each computing job is characterized by two key attributes: the number of servers required and the job’s running duration (in hours). We introduce two finite sets: $\mathcal{K} = \{1, 2, 4, 8, 16\}$ denotes the set of possible server requirements, and $\mathcal{L} = \{1, 2, 4, ..., 168\}$ denotes the set of possible running durations. We define $\mathcal{KL} = \mathcal{K} \times \mathcal{L}$ as the set of all possible job-type combinations (server demand, running time). Based on this structure, we aggregate the jobs by grouping together those with identical resource demands ($k \in \mathcal{K}$) and identical running durations ($l \in \mathcal{L}$). We define $T_h$ as the optimization time for each round. For each time $r \in T$, the optimization time domain is represented as $T_r = \{r, r + 1, ..., r + T_h - 1\}$. We further assume that all jobs are independent and that their execution order does not affect overall system performance.

\subsection{Data Center Job Scheduling Model}
Our goal is to maintain a high level of load completion while minimizing the operational cost. To ensure the fundamental operational requirements of the DC, we include a job-completion term in the objective function together with penalties on peak demand and electricity usage. Our aim is to complete incoming jobs as promptly as possible while minimizing both energy expenses and peak electricity demand, thereby reducing stress on the grid’s transmission network. We first consider:
\begin{equation}\label{active}
\begin{aligned}
m(t) &= \sum_{(k,l) \in \mathcal{K} \mathcal{L}} \sum_{t' = st(t,l)}^{t} k \, n_{kl}(t') + \sum_{l = t-r+1}^{\bar L - 1} u_l(t) \\
&\quad - \sum_{(t_b,k,l) \in V(t)} k \, v^{t_b}_{kl}(r), 
\quad \forall t \in \mathcal{T}_r,
\end{aligned}
\end{equation}

\begin{equation}\label{capBound}
	\begin{aligned}
	m(t) \leq \bar I_t, 
    && \forall t \in \mathcal{T}_r.
	\end{aligned}
\end{equation}
where $s_t(t, l) = \max\{r, t - l + 1\}$ represents the earliest start time of a job of length $l$ at time $t$. $\bar L = \max \{L\}$ and $\bar I_t = \max \{I_t\}$. $V(t) := \{(t_b, k, l) \mid t_b + l - t > 0,\ (k,l) \in \mathcal{KL}\}$ denotes the set of all tuples $(t_b, k, l) $such that jobs of type ($k, l)$ started at time $t_b$ have not yet finished by time $t$. The decision variables $m(t)\in\mathbb{Z}_{+}$ and $n_{kl}(t')\in\mathbb{Z}_{+}$ represent the number of active servers at time $t$ and the number of $(k,l)$ type jobs started at time $t'$, respectively. The decision variable $v^{t_b}_{kl}(r)\in\mathbb{Z}_{+}$ denotes the number of \((k,l)\) jobs that are canceled at the beginning of time \(r\) and were originally started at time \(t_b\). The state variable $u_l(t)\in\mathbb{Z}_{+}$ represents the number of servers occupied by jobs with a remaining time of $l$ hours at time $t$. Equation \eqref{active} defines the number of active servers, and Equation \eqref{capBound}  represents the server capacity constraint.

Then, we introduce the basic constraints of job scheduling. We define $N^e_{kl}(t')\in\mathbb{Z}_{+}$ as the expected number of $(k,l)$ type jobs at time $t'$. The state variables $\tilde{n}_{kl}(r)\in\mathbb{Z}_{+}$ and $\hat n_{kl}^{t_b}(r)\in\mathbb{Z}_{+}$ are the number of $(k,l)$ type jobs in the waiting queue at time $r$, and the number of $(k,l)$ type jobs that started at time $t_b$ and are still running at time $r$. Equations \eqref{clearance}–\eqref{cancelbound} define workload execution limits, backlog constraints, and job-cancellation feasibility.
\begin{equation}\label{clearance}
	\begin{aligned}
	\sum_{t'=r} n_{kl}(t') \leq \sum_{t'=r}N^e_{kl}(t') + \tilde{n}_{kl}(r), 
    && \forall (k,l) \in \mathcal{KL},
	\end{aligned}
\end{equation}
\begin{equation}\label{minjobclear}
    \begin{aligned}
     \sum_{t'=r}^{r+T_h -1} n_{kl}(t') \leq \sum_{t'=r}^{r+[T_h/2] -1} N^{e}_{kl}(t') + \tilde n_{kl}(r), && \forall (k,l) \in \mathcal{KL},      
    \end{aligned}
\end{equation}
\begin{equation}\label{cancelbound}
    \begin{aligned}
    v_{kl}^{t_b}(r) \leq \hat n_{kl}^{t_b}(r), && \forall t \in \mathcal{T}_r, \forall (k,l) \in \mathcal{KL}.     
    \end{aligned}
\end{equation}

We now discuss the constraints regarding power consumption for the DC, we have:
\begin{equation}\label{gridpower1}
    \begin{aligned}
    P(m(t)) = g(t) = (P_p - P_i) \frac{m(t)}{I} + P_i,
    && \forall t \in \mathcal{T}_r,  
    \end{aligned}
\end{equation}

\begin{equation}\label{cost}
    \begin{aligned}
    \psi (t) = g(t) \cdot p(t), && \forall t \in \mathcal{T}_r, 
    \end{aligned}
\end{equation}

\begin{equation}\label{peakpower}
    \begin{aligned}
     G(r) \geq g(t), \quad \forall t \in \mathcal{T}_r.
     \end{aligned}
\end{equation}
where $P(m(t))\in\mathbb{R}_{+}$ and $g(t)\in\mathbb{R}_{+}$ represent the electricity consumption and electricity purchased of DC from the grid at time $t$. $p(t)\in\mathbb{R}$ is the electricity price at time $t$, while $\psi (t)\in\mathbb{R}$ represents the DC usage charge. $G(r)\in\mathbb{R}_{+}$ represents the peak power at the current interval $r$. $P_p\in\mathbb{R}_{+}$ and $P_i\in\mathbb{R}_{+}$ are defined as the peak and idle power of DC. 

For a detailed scheduling model, please refer to \cite{zhang2022exploring}.

\subsection{Energy Management with On-site Solar Power}
After completing the DC job scheduling model, we considered the impact of adding on-site solar on the DC. We introduced the decision variable $s(t)\in\mathbb{R}_{+}$ and parameter $\xi(t)\in\mathbb{R}_{+}$, which represent the solar power consumed by the DC at time $t$ and the maximum energy that the on-site solar system can provide. Since the on-site solar system can be started and stopped within minutes \cite{loutan2017demonstration}, meaning it can freely discard power, we have:
\begin{equation}\label{solar31}
    \begin{aligned}
	0 \leq s(t) \leq \xi(t),
    && \forall t \in \mathcal{T}_r.
    \end{aligned}
\end{equation}

Consequently, in contrast to Equation \eqref{gridpower1}, the electricity procured from the power grid by the system will be given by:
\begin{equation}\label{gridpower3}
    \begin{aligned}
    g(t) = (P_p - P_i) \frac{m(t)}{I} + P_i - s(t),
    && \forall t \in \mathcal{T}_r.
    \end{aligned}
\end{equation}

That is, the amount of electricity purchased by the DC from the grid at time $t$ equals the total power consumption minus the solar power utilized. If we want to consider that DC prioritizes the consumption of its own generated solar power, we need to set a constraint that the solar power $s(t)$ used by DC at time $t$ is the minimum between the total energy consumption $P(m(t))$ and the maximum solar power provided $\xi(t)$. We have:
\begin{equation}\label{solar0}
\begin{aligned}
s(t) = \min\{P(m(t)), \ \xi (t)\},
&& \forall t \in \mathcal{T}_r. 
\end{aligned}
\end{equation}

To linearize Equation \eqref{solar0}, we introduce a big-M formulation. Specifically, if $P(m(t)) \geq \xi(t)$, then $y(t)=1$ and $s(t)=\xi(t)$. If $P(m(t)) < \xi(t)$, then $y(t)=0$ and $s(t)=P(m(t))$, where $M$ is selected to be sufficiently large. From this formulation, we get:
\begin{equation}\label{solar3}
	\begin{aligned}
	s(t) \geq \xi(t) - M \cdot (1-y(t)),
    && \forall t \in \mathcal{T}_r,
	\end{aligned}
\end{equation}
\begin{equation}\label{solar4}
	\begin{aligned}
	s(t) \geq P(m(t)) - M \cdot y(t),
    && \forall t \in \mathcal{T}_r,
	\end{aligned}
\end{equation}
\begin{equation}\label{solar2}
	\begin{aligned}
	s(t) \leq P(m(t)), 
    && \forall t \in \mathcal{T}_r,
	\end{aligned}
\end{equation}
\begin{equation}\label{solar5}
	\begin{aligned}
	y(t) \in \{0,1\}.
	\end{aligned}
\end{equation}

We further consider that if a feedback mechanism is introduced to sell excess solar power generated by the DC back to the grid, then we have:
\begin{equation}\label{feedpower}
    \begin{aligned}
    f(t) &= \xi(t) - s(t), && \forall t \in \mathcal{T}_r,
    \end{aligned}  
\end{equation}
\begin{equation}\label{Case4cost}
    \begin{aligned}
    \psi (t) = g(t)\cdot p(t) - f(t)\cdot x(t),
    && \forall t \in \mathcal{T}_r.
    \end{aligned}
\end{equation}
where the decision variable $f(t)\in\mathbb{R}_{+}$ denotes the amount of power sold to the grid at time $t$, and parameter $x(t)\in\mathbb{R}_{+}$ denotes the grid purchase price. 

\subsection{Battery Participation in Flexible Scheduling}
After incorporating solar power, we aim to further enhance the flexibility of the DC scheduling system by introducing the battery. Consequently, the electricity supply of DC will have three sources: solar generation, battery storage, and the power grid. The energy generated by the solar system can be flexibly allocated either directly to the DC, stored in the battery, or sold back to the grid. We consider:
\begin{equation}\label{battery1}
    \begin{aligned}
	Q(t+1) = Q(t) + \eta_{\text{in}} \cdot b_{\text{in}}(t) - \frac{1}{\eta_{\text{out}}} \cdot b_{\text{out}}(t),
    && \forall t \in \mathcal{T}_r,
    \end{aligned}
\end{equation}
\begin{equation}\label{battery2}
    \begin{aligned}
    0 \leq Q(t) \leq \bar Q(t), 
    && \forall t \in \mathcal{T}_r,
    \end{aligned}
\end{equation}
\begin{equation}\label{battery3}
    \begin{aligned}
     0 \leq b_{\text{in}}(t) \leq \bar b_{\text{in}}, \quad \forall t \in \mathcal{T}_r,
     \end{aligned}
\end{equation}
\begin{equation}\label{battery4}
    \begin{aligned}
     0 \leq b_{\text{out}}(t) \leq \bar b_{\text{out}}, \quad \forall t \in \mathcal{T}_r.
     \end{aligned}
\end{equation}   
where the state variable $Q(t)\in\mathbb{R}_{+}$ denotes the energy level of the battery at time $t$. The decision variables $b_{\text{in}}(t)\in\mathbb{R}_{+}$ and $b_{\text{out}}(t)\in\mathbb{R}_{+}$ represent the amount of energy charged into and discharged from the battery at time $t$. $\bar Q(t)\in\mathbb{R}_{+}$, $\bar b_{\text{in}}\in\mathbb{R}_{+}$, and $\bar b_{\text{out}}\in\mathbb{R}_{+}$ denote the maximum energy capacity of the battery, the maximum charging power, and the maximum discharging power. The parameters $\eta_{\text{in}}\in\mathbb{R}_{+}$ and $\eta_{\text{out}}\in\mathbb{R}_{+}$ are the charging and discharging efficiency coefficients. Equation \eqref{battery1} illustrates the state update of the battery. To prevent the battery from charging and discharging simultaneously at the same time, we have:
\begin{equation}\label{battery5}
    \begin{aligned}
     b_{\text{in}}(t) \cdot b_{\text{out}}(t) = 0, \quad \forall t \in \mathcal{T}_r.
     \end{aligned}
\end{equation}   

To linearize it while preserving the model's solvability, we introduce a binary variable \( z(t) \in \{0,1\} \) for each time period \( t \), which indicates the operating mode of the battery:

\begin{equation}\label{battery6}
    \begin{aligned}
     b_{\text{in}}(t) \leq \bar{b}_{\text{in}} \cdot z(t), \quad \forall t \in \mathcal{T}_r,
     \end{aligned}
\end{equation}   
\begin{equation}\label{battery7}
    \begin{aligned}
     b_{\text{out}}(t) \leq \bar{b}_{\text{out}} \cdot (1 - z(t)), \quad \forall t \in \mathcal{T}_r.
     \end{aligned}
\end{equation}   

Meanwhile, the amount of electricity purchased by the DC from the grid at time $t$ becomes the total power consumption minus the solar power used by the DC and the energy discharged from the battery. Similarly, the amount of electricity sold back to the grid by the DC at time $t$ equals the total solar generation minus the solar power consumed by the DC and the energy charged into the battery. We have:
\begin{equation}\label{power51}
    \begin{aligned}
    g(t) = (P_p - P_i) \frac{m(t)}{I} + P_i- s(t) - b_{\text{out}}(t) ,
    && \forall t \in \mathcal{T}_r,
    \end{aligned}
\end{equation}
\begin{equation}\label{power52}
    \begin{aligned}
    f(t) = \xi (t) - s(t) - b_{\text{in}}(t),
    && \forall t \in \mathcal{T}_r.
    \end{aligned}
\end{equation}

\subsection{Optimization Problem and Objective Function}

This study aims to develop an integrated computational framework that enables joint optimization of DC job scheduling and energy management under dynamic conditions including volatile electricity prices, fluctuating solar generation, and time-varying job arrivals. We first consider incentives for proactive job scheduling and guarantees on job completion, leading to the following equation:
\begin{equation}\label{UTIL}
\begin{split}
\theta(r) &= \sum_{t \in \mathcal{T}_r} \sum_{(k,l) \in \mathcal{KL}} \big((r + T_h) k l - t\big) n_{kl}(t) \\
&\quad - \sum_{(k,l) \in \mathcal{KL}} \sum_{t_b = T_{rl}}^{r-1} \big((r + T_h) k l - t_b\big) v^{t_b}_{kl}(r)
\end{split}
\end{equation}
where $T_{rl}$ represents the earliest start time of the job of $(k,l)$ type. The decision variable $n_{kl}(t)\in\mathbb{Z}_{+}$ represents the job of $(k,l)$ type starting at time $t$. The first term on the right-hand side of the equation shows the reward for proactive job scheduling, with $(r + T_h) k l - t$ as the reward weight, meaning the earlier the job is scheduled, the greater the reward. The second term shows the penalty for job cancellation, with $(r + T_h) k l - t_b$ as the penalty weight, meaning the later the job is cancelled, the greater the penalty. Meanwhile, we introduce $\lambda_p\in\mathbb{R}_{+}$  as the penalty coefficients for peak power demand, yielding the following RHO model at time intervals $r$ over the decision horizon $\mathcal T_r$:
 \begin{equation}\label{case1Obj}
	\begin{aligned}
	\max_{m, n, u, v, s} \; \theta(r) - \lambda_p \cdot G(r) -\sum_{t\in\mathcal T_r} \psi(t)
	\end{aligned}
\end{equation}
\begin{equation}\label{others}
    \textrm{s.t.}\;
    \begin{aligned}
    \eqref{active}\text{--}\eqref{UTIL}.
    \end{aligned}
 \end{equation}
\section{Case Studies}
We select a CAISO node (35.56°N, 118.96°W) with high annual solar irradiance (234.58~Wh/m$^2$, which is higher than the state average 216.26~Wh/m$^2$) and large price volatility (mean 57.77 USD/MW, variance 2337.18). Five operational strategies were simulated and compared:
\begin{enumerate}[label=(\roman*)]
    \item \textbf{Grid-only operation with peak and usage penalties (Case I):}
    the DC procures all electricity from the grid. The scheduling decisions respond only to usage charges (time-varying electricity prices) and peak-demand penalties; no on-site solar, feed-in, or storage is available.
    \item \textbf{Prioritized solar consumption (Case II):}
    on-site solar is used under a strict self-consumption priority. Whenever solar is available, the DC is forced to consume solar up to its instantaneous power demand, and any residual demand is met by grid procurement; solar export is not allowed.
    \item \textbf{Flexible solar--grid co-optimization (Case III):}
    solar utilization is a decision variable. The model jointly optimizes how much solar to use and how much grid electricity to procure at each time step, allowing solar curtailment when economically optimal; solar export is not allowed.
    \item \textbf{Surplus solar feed-in to the grid (prosumer mode) (Case IV):}
    surplus solar can be exported to the grid under a feed-in (LMP) price. The model jointly optimizes solar self-consumption, solar export, and grid procurement, enabling the DC to operate as a prosumer.
    \item \textbf{Battery participation in flexible scheduling (Case V):}
    a battery is added to provide temporal shifting and additional flexibility. Solar can be allocated among direct DC use, battery charging, and grid export, while the battery can be discharged to reduce grid purchases and manage peak demand.
\end{enumerate}

\begin{figure*}[!t]
	\center{\includegraphics[width=0.8\textwidth]{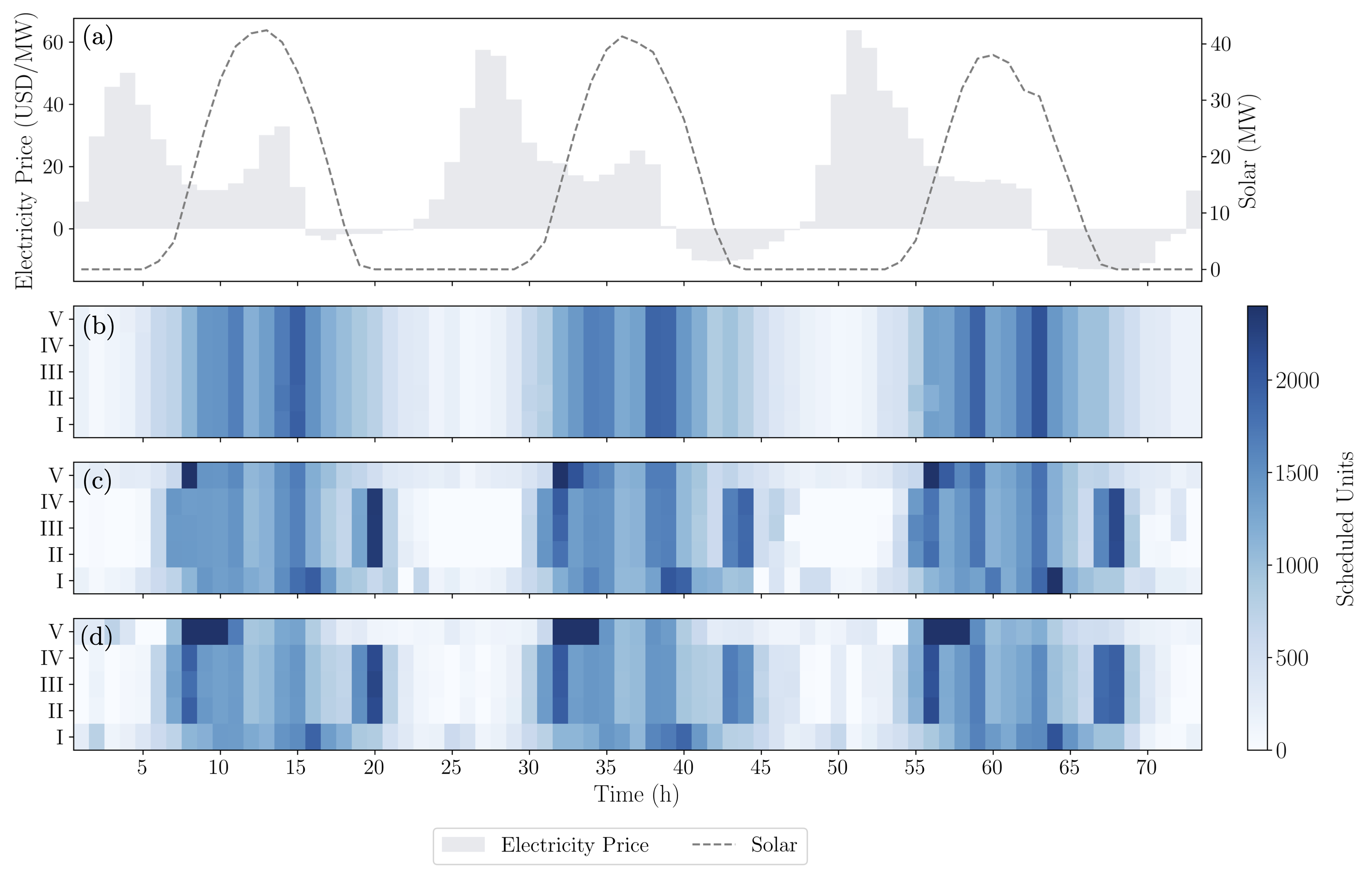}}
	\caption{Job scheduling under different peak demand penalties. Darker shades indicate higher scheduling intensity. (a) shows the changes in electricity price and solar power over time. (b)-(d) represent job scheduling when \(\lambda_p = 0\), 15, and 50. The DC schedules jobs based solely on electricity price when \(\lambda_p = 0\). As \(\lambda_p\) increases, the impact of solar power begins to manifest.}
	\label{scheduled}
    \vspace{1em}
    \center{\includegraphics[width=0.85\textwidth]{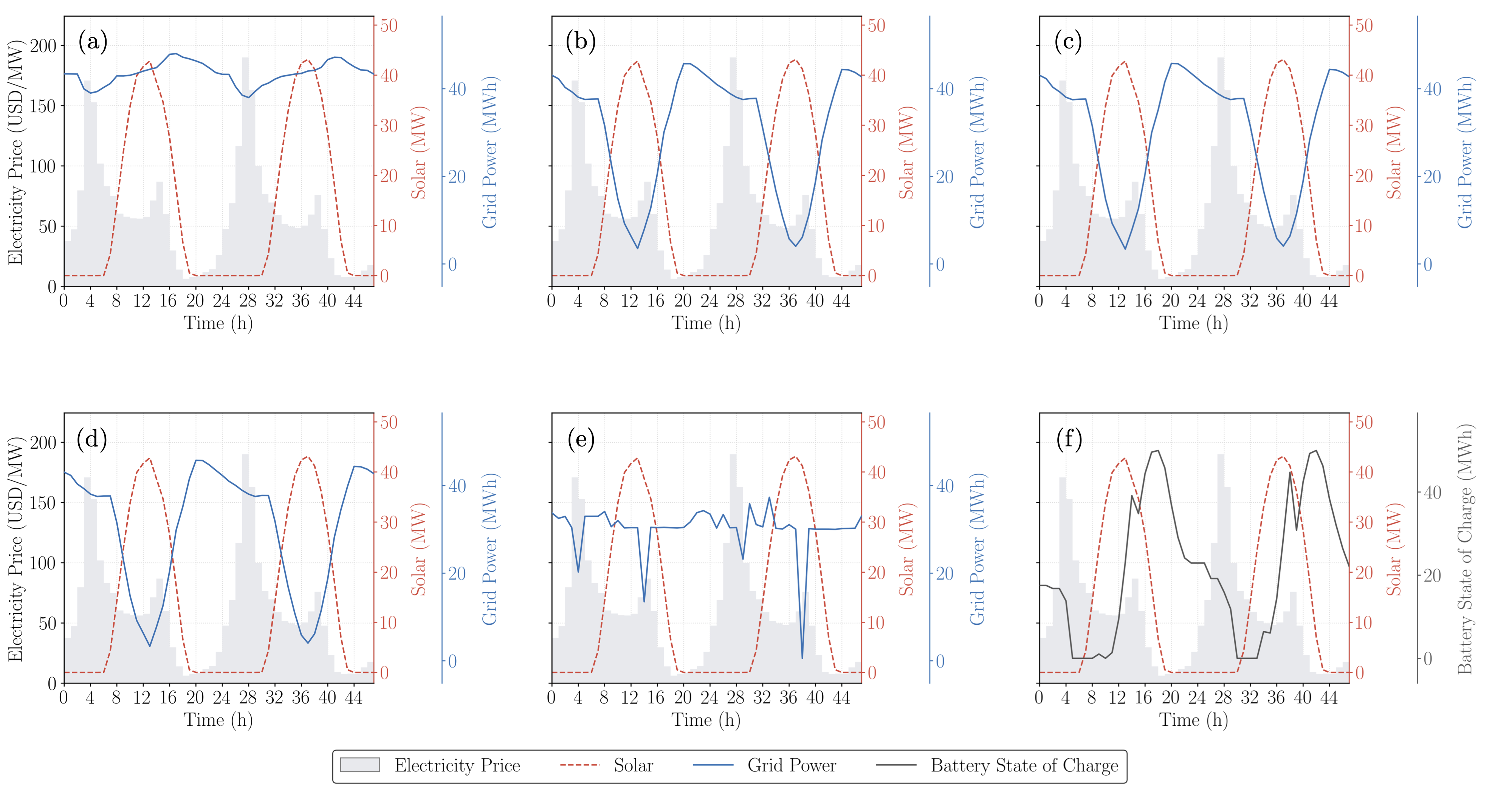}}
	\caption{Electricity purchasing pattern of the DC over time. Gray bars are electricity prices, and red lines are solar generation. In (a)–(e), blue lines denote the amount of grid electricity purchased in Cases I–V, respectively. In (f), the dark gray line represents the battery state of charge in Case V. The addition of on-site solar will not make the DC unresponsive to electricity prices, and battery storage provides greater flexibility.}
	\label{gridpower}

\end{figure*}

\begin{figure*}[!t]
\center{\includegraphics[width=0.85\textwidth]{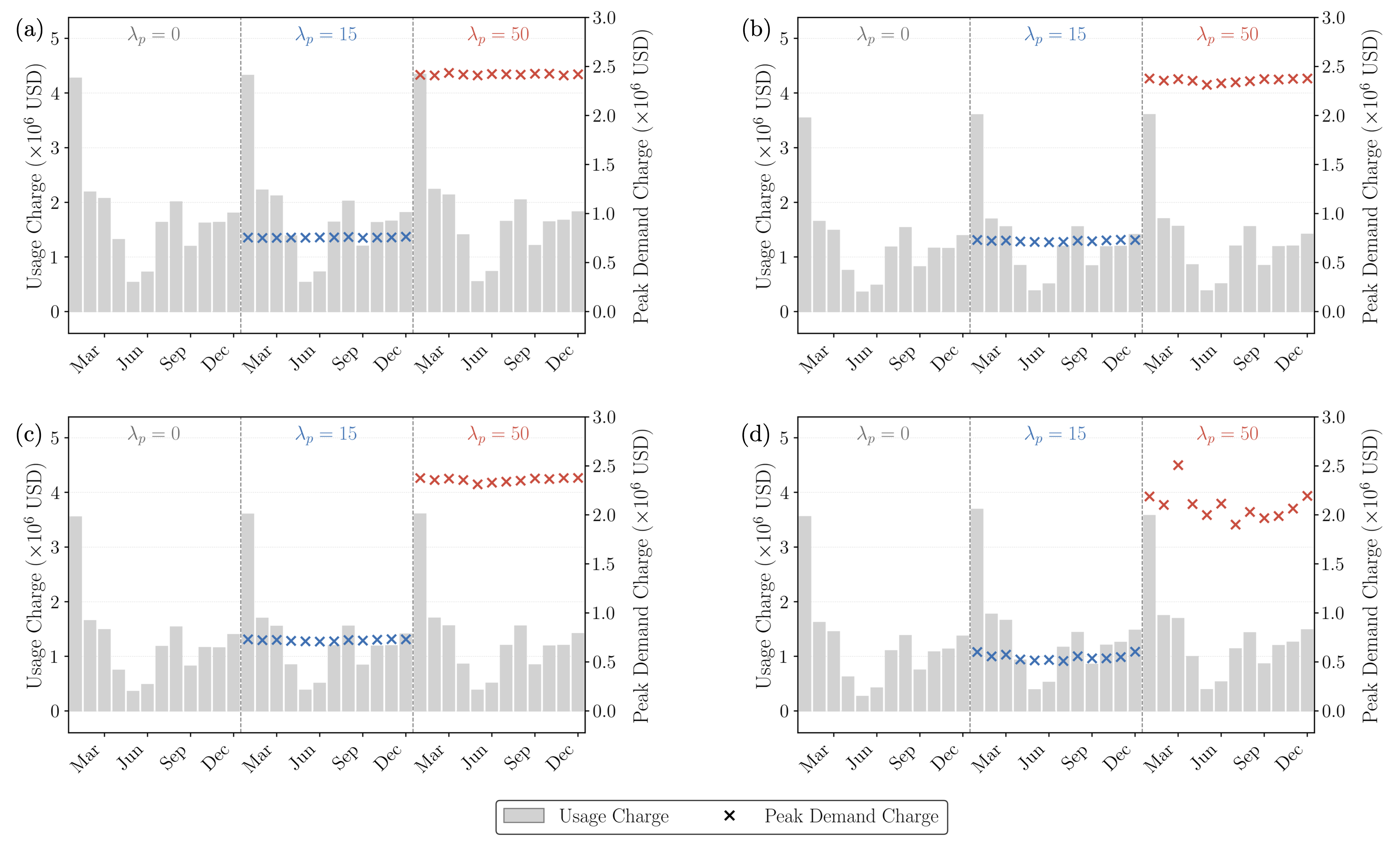}}
	\caption{Electricity cost of the DC under different \(\lambda_p\). The gray bars represent the usage charge, and the scatter points denote the peak demand charge. Each subfigure is divided into three sections corresponding to different values of \(\lambda_p\). (a)–(d) illustrate the cost distribution for Case I, III, IV, and V, respectively. Since the cost variations in Case II and Case III are nearly identical, Case III is used as the representative in (b). On-site solar can effectively reduce usage charges, but it has little effect on peak demand charges.}
	\label{ecost}
  
\end{figure*}

\begin{figure}[!htb]
	\center{\includegraphics[width=0.99\columnwidth]{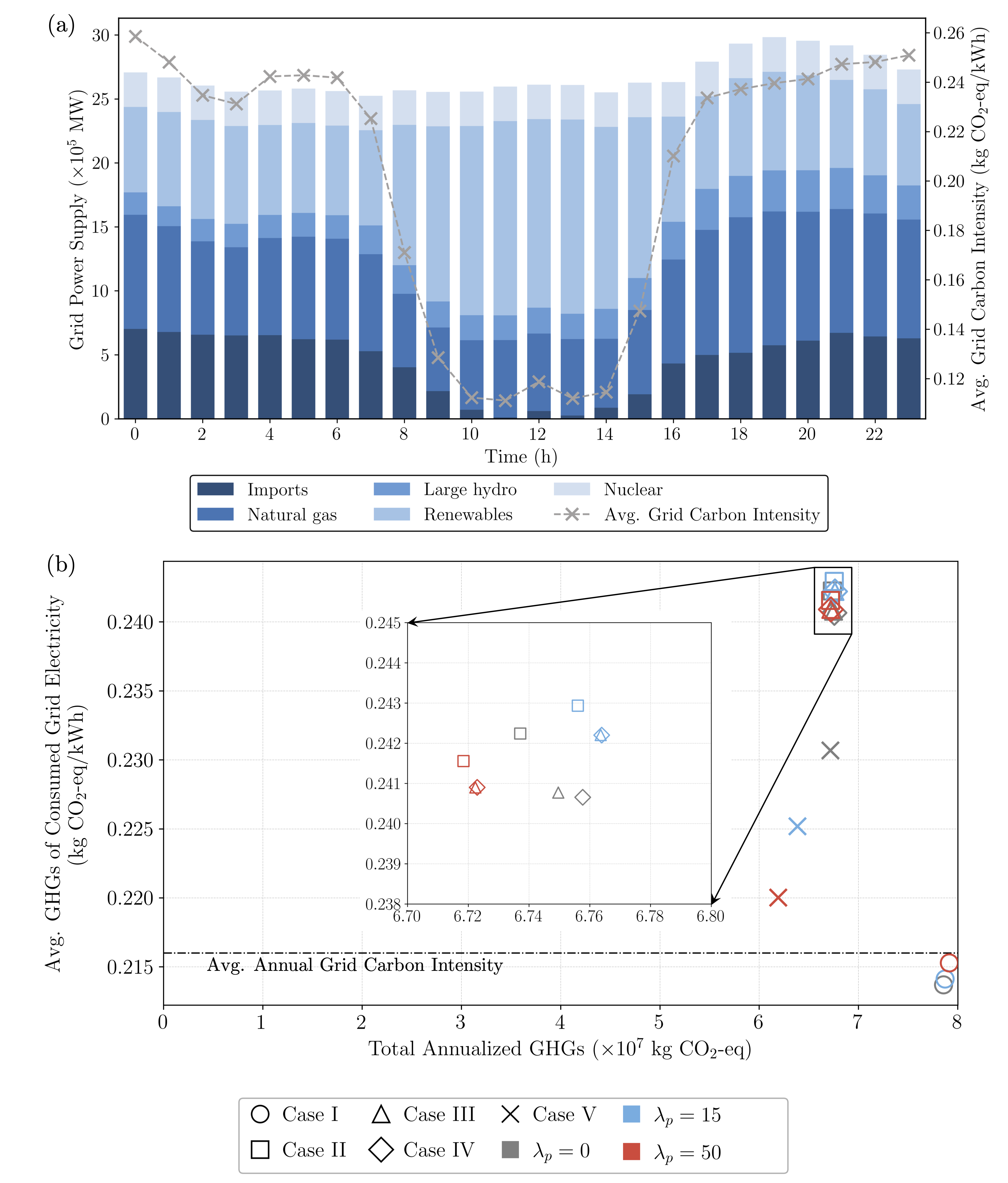}}
	\caption{Carbon emissions analysis. (a) shows the sources of power generation for CAISO and the average carbon intensity of the power grid over a day. (b) shows the trade-off between total annual carbon emissions and the average carbon intensity of grid electricity consumed by DC.}
	\label{carb}
\end{figure}
In our case study, we consider a large-scale DC with a capacity of 20,000 servers, a peak power of 100 MW, and an idle power of 30 MW. We set the solar capacity to 50 MW \cite{shehabi20242024, luthander2019graphical} with a battery capacity to solar capacity ratio of 1:1 \cite{BatterySkills_SolarBatteryCapacity}. We run a one-year closed-loop simulation (168 h horizon, 1 h resolution) in Julia 1.10.5~\cite{bezanson2017julia} with JuMP 1.23.2~\cite{dunning2017jump} and Gurobi 11.0.2 ~\cite{Gurobi_DecisionIntelligence_Leader}.

The values of \(\lambda_p\) are selected based on the current electricity market conditions. The CAISO imposes relatively high peak demand penalties (approximately 50\,USD/kWh), while regions such as Texas (ERCOT) and parts of the Midwest do not apply such penalties. Accordingly, we consider \(\lambda_p = 0\), 15 and 50, where \(\lambda_p = 15\) represents an national average penalty level \cite{mclaren2017identifying}. We find increasing \(\lambda_p\) introduces a small number of job cancellations and slightly longer delays because the DC tends to defer workloads to avoid peak demand charges. Cases with on-site solar exhibit marginally lower scheduling flexibility, whereas prosumer operation partially alleviates this effect. Despite these trade-offs, service quality remains high, with job completion rates consistently above 99.9\% in all simulations.

Fig.~\ref{scheduled} shows the job scheduling intensity of three arbitrary days with different level of peak demand penalty. Due to the periodicity of solar generation, electricity prices, and job input, the scheduling results also exhibit clear periodic patterns. When \(\lambda_p = 0\), since solar power is insufficient to support power consumption most of the time, the scheduling of jobs for all cases strictly follows the power response to reduce the electricity bill. When imposing a peak demand penalty (\(\lambda_p \geq 0 \)), job scheduling remains flexible and responds to both electricity prices and solar availability. The Case~I mode smooths the job deployment evenly and results in the lowest peak penalty. Interestingly, in Cases II to IV, the DC not only responds to the volatility of the solar power but also schedules a large number of jobs to the night. This is because during the daytime high-price period, the solar can only offset part of the demand. Over-scheduling would result in high costs, so the system chooses to sacrifice peak demand charge to reduce overall costs. The introduction of the battery will cause jobs to be scheduled more frequently during the daytime when electricity prices are low and solar power is plentiful. The storage characteristics allow it to avoid periods of high electricity prices and make full use of solar power, while avoiding excessive peak penalties and achieving minimal cost.

In Fig.~\ref{gridpower}, the electricity purchasing pattern in Case I is clearly anti-correlated with the price trend: grid usage falls when prices rise and increases when prices drop. In Cases II and III, the on-site solar power offsets a large portion of the grid's electricity during the daytime, which dents the electricity purchasing curve. Cases~II and~III not only consume substantially less grid electricity than Case~I but also exhibit very similar purchasing patterns. This suggests that regardless of whether a priority constraint on solar utilization is imposed, the optimized decisions naturally favor available solar generation in most situations. Meanwhile, the overall consumption trend in these cases remains negatively correlated with electricity prices. As solar capacity does not exceed that of DC and the solar power is absorbed by DC itself, very little electricity is exported to the grid, so the electricity purchasing pattern in Case IV is similar to that in Cases II and III. For Case~V, the introduction of the battery further enhances flexibility in grid interaction. The purchasing pattern becomes more volatile, reflecting the DC’s improved capability to adjust and identify the optimal energy mix. The inclusion of storage allows the DC to store energy when solar generation is abundant and electricity prices are low, and to supply power when solar output is insufficient and prices are high.

In Fig.~\ref{ecost}, we compare the usage charge and peak demand penalty across operational modes. In the first four cases, the peak demand charges remain perfectly flat throughout the year when penalties are enforced ($\lambda_p = 15$ and $50$), as the facility hits the exact same maximum instantaneous power draw every single month. The introduction of on-site solar reduces monthly usage charges, but it does not lower peak demand by too much even with a high peak demand penalty factor as $\lambda_p = 50$. The system trade off peak demand and usage charges in order to achieve the lowest electricity cost. Being flexible with electricity prices and ramping up the system when the price is low provides a large incentive, even if it causes a slightly higher peak demand. Case V demonstrates that adding a battery fundamentally changes this dynamic. Energy storage provides dispatch flexibility, allowing the system to actively discharge during the specific hours of maximum net load. This targeted peak shaving successfully suppresses the monthly peak demand charges, with the most savings achieved during the summer months when abundant solar power reliably recharges the battery to offset cooling and IT peaks.

\begin{figure*}[t]
	\center{\includegraphics[width=0.85\textwidth]{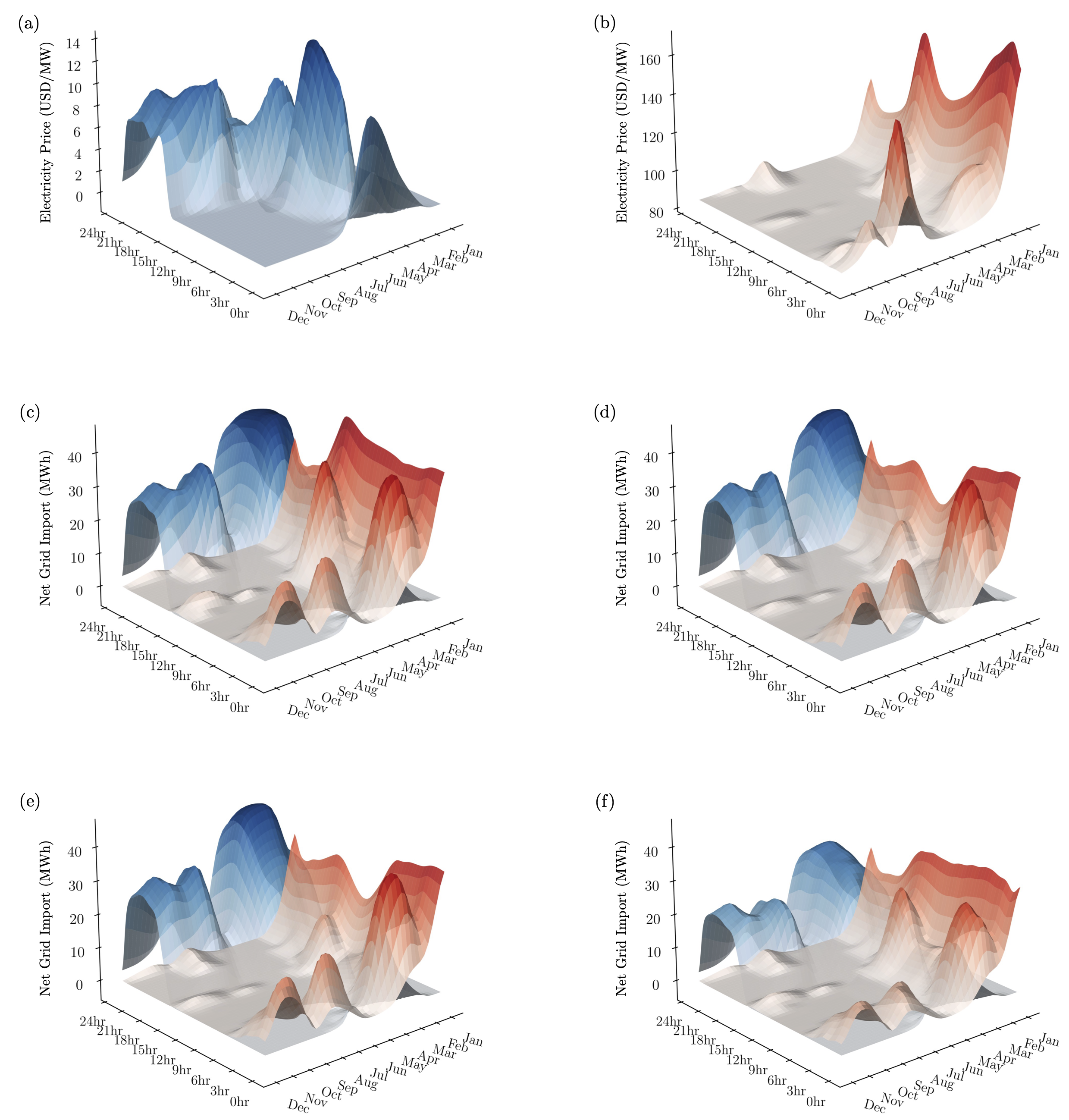}}
	\caption{Electricity transactions during peak and off-peak hours. (a) and (b) show the distribution of electricity prices during off-peak and peak hours. (c)–(f) show the net grid import from the grid in Cases I, III (Cases II and III have the identical patterns, and only one plot is presented here), IV, and V during peak hours (red) and off-peak hours (blue). Incorporating an on-site solar facility, DCs discharge less power from the grid during peak hours while injecting a larger amount of electricity into the grid during off-peak hours.}
	\label{flex}
\end{figure*}

In Fig.~\ref{carb}, the carbon intensity of grid electricity exhibits significant intra-day variation driven by the real-time generation mix. This overall intensity is determined by the fractional contribution of each energy source and its specific emission factor. A diurnal pattern is the marked decrease in carbon intensity during midday hours, which directly corresponds to peak solar generation on the grid. To decipher whether the on-site solar affects the DC's ability to absorb low-carbon electricity from the grid, we present the total GHGs of the DC, accounting for electricity-related carbon emissions, the annualized life-cycle emissions of the solar PV and battery systems, and the average carbon intensity of the consumed grid electricity. The baseline average carbon intensity of electricity generation in CAISO in 2023 is \(0.216\,\mathrm{kg\ CO_{2}\text{-}eq}\) per kWh~\cite{CAISO_TodaysOutlook}. A comparison reveals that in Case~I, active job scheduling enables the DC to consume a relatively high share of low-carbon grid electricity. When on-site solar is introduced, the DC’s overall energy mix becomes cleaner; however, the remaining electricity purchased from the grid tends to have a higher average carbon intensity. This occurs because periods of abundant on-site solar generation in the CAISO region directly coincide with a high share of utility-scale solar on the broader grid. Consequently, strictly prioritizing on-site solar self-consumption during these hours displaces the DC's opportunity to absorb available green electricity from the grid. The prosumer mode slighty mitigate this unintended effect. While adding battery storage lowers the DC's overall emissions without compromising its ability to utilize low-carbon grid power too much, especially when there is a high penalty on peak demand. Although the annualized life cycle emissions of the on-site solar system amount to \(1.65 \times 10^{6}\,\mathrm{kg\ CO_{2}\text{-}eq}\) (estimated based on data from ~\cite{danelli2024environmental}), the total DC emissions remain lower than those of Case~I after solar integration. The battery system generates \(2.33 \times 10^{5}\,\mathrm{kg\ CO_{2}\text{-}eq}\) over its annualized life cycle (estimated based on data from~\cite{IEA_EVBatterySupplyChainSustainability_2024, cole2025cost}); therefore, the annual total emissions in Case~V are also at a low level. In summary, the introduction of on-site solar can lower the total annual GHGs of the DC by roughly 10-15\%, but still constrains the DC's ability to absorb low-carbon electricity from the grid, even with a flexible operating mode such as prosumer. The reduction of total annual GHGs comes at a cost of more than a 10\% increase in the average carbon intensity of the consumed grid electricity.

To reveal the impact of an on-site solar facility on the DC's ability to serve as a temporally flexible asset for the grid, we set the average DAM price of CAISO in 2023 for the studied region as the 50\% reference line. When the DAM price falls within the range of 25\%–75\% of this average value, the electricity market is defined as operating under normal conditions, while all remaining periods are classified as abnormal~\cite{luoma2014forecast}. The abnormal periods are further divided into two categories: when the electricity price is below 25\%, the grid is considered electricity sufficient (off-peak hour), whereas when the price exceeds 75\%, the grid is regarded as electricity insufficient (peak hour)~\cite{kohansal2017data}. Fig.~\ref{flex} presents the electricity prices and net grid import (purchased electricity minus sold electricity) of the DC during peak and off-peak hours. Figs.~\ref{flex}(a) and \ref{flex}(b) show that off-peak hours mainly occur at night as social electricity consumption is low and during the summer as renewable energy is abundant. Peak hours mainly occur in the early morning and during winter. Particularly, CAISO electricity prices were unusually high in January 2023, with a monthly average price of 138.56\,USD/MW~\cite{CAISO_TodaysOutlookPrices_2025}. For Case~I, the DC has the highest frequency of charging from the grid during peak hours in January. Abnormal prices make it unavoidable for the DC to purchase high-priced electricity. The integration of on-site solar effectively reduces electricity purchases during peak hours from 9 AM to 4 PM, when solar is abundant, thereby alleviating grid stress. The battery storage decreases and smooths peak-hour grid electricity usage. Instead of only backing up the grid during the daytime when solar is available (from 9 AM to 4 PM), the system also alleviates grid stress in the early morning and at the evening, since the battery adds another layer of temporal flexibility to the system. In Fig.~\ref{flex}(e) –~\ref{flex}(f), the blue shade represents the frequency of net grid import during off-peak hours. One can observe that DC tends to discharge a slightly lower amount of electricity from the grid during off-peak hours when integrating on-site solar. This is because more computing jobs were completed in the daytime when solar is on (since solar is free), and fewer computing jobs were carried out during the nighttime even though the electricity price is cheap. Especially when integrated with a battery, this pattern becomes more obvious: the battery tends to store the excess solar power, and the profit-oriented schedule model prioritizes the use of stored solar power instead of off-peak electricity from the grid. In summary, the addition of on-site solar alleviates the grid stress in daytime (peak shaving) but also constrains DC from absorbing off-peak electricity since the DC is soaking into the on-site solar.

\section{Conclusion}  

This paper developed a receding-horizon optimization framework to coordinate data-center workload scheduling, grid electricity procurement, on-site solar utilization, solar export, and battery operation. Using a year-long CAISO case study, we show that behind-the-meter solar reduces electricity costs and total facility-level greenhouse gas emissions, but has limited impact on peak demand charges. More importantly, solar self-consumption can increase the average carbon intensity of the remaining grid electricity consumed by the data center, because local solar generation often coincides with low-carbon periods in the bulk grid. This reveals a trade-off between local renewable self-sufficiency and system-level decarbonization. Prosumer operation and battery storage partially mitigate this trade-off by enabling solar export, temporal energy shifting, and peak shaving. These results highlight the need to coordinate local renewable utilization with grid-level price and carbon signals when designing sustainable data centers.


\bibliography{ref}
\bibliographystyle{IEEEtran}

\vfill

\end{document}